\documentclass[12pt,leqno]{amsart}
\usepackage{amsfonts}                   
\usepackage{latexsym}                   
\usepackage{amssymb}
\usepackage[lite]{amsrefs}
\usepackage{graphics}
\usepackage{graphics}
\usepackage[lite]{amsrefs}
\usepackage{amsfonts}
\usepackage{amssymb}
\usepackage{xypic,amsmath,amssymb}
\usepackage{amsmath}
\usepackage[dvips]{graphicx}                   
\newtheorem{teorema}{Theorem}[section] 
\newtheorem{lemma}[teorema]{Lemma}
\newtheorem{propos}[teorema]{Proposition}
\newtheorem{corol}[teorema]{Corollary}
\newtheorem{ex}{Example}[section]
\newtheorem{rem}{Remark}[section]
\newtheorem{defin}[teorema]{Definition}

\def\defin{\par\ifdim\lastskip<\smallskipamount\removelastskip
  \smallskip\fi\noindent{\bf\ignorespaces Definition\unskip:\enspace}\rm
\ignorespaces}

\def\bit{\begin{itemize}}
\def\eit{\end{itemize}}
\def\be{\begin{equation}}
\def\ee{\end{equation}}
\def\beq{\begin{eqnarray}}
\def\eeq{\end{eqnarray}}

\def\ba{\begin{array}}
\def\ea{\end{array}}
\def\bt{\begin{teorema}}
\def\et{\end{teorema}}
\def\bp{\begin{propos}}
\def\ep{\end{propos}}
\def\bl{\begin{lemma}}
\def\el{\end{lemma}}
\def\bc{\begin{corol}}
\def\ec{\end{corol}}
\def\br{\begin{rem}\rm}
\def\er{\end{rem}}
\def\bex{\begin{ex}\rm}
\def\eex{\end{ex}}
\def\bd{\begin{defin}}
\def\ed{\end{defin}}
\def\demo{\par\noindent{\bf Proof.\ }}
\def\enddemo{\ $\Box$\par\vskip.6truecm}

     \def\nin{\noindent}

\def\ES{\varnothing}  \def\IN{\infty}  

\def\R{{\mathbb {R}}}   \def\a {\alpha} \def\b {\beta}\def\g{\gamma}
\def\N{{\mathbb N}}      \def\e{\varepsilon}
\def\C{{\mathbb C}}      \def\l{\lambda}
					\def\o{\omega}\def\p{\partial}\def\r{\varrho}
\def\Z{{\mathbb Z}}					\def\s{\sigma}\def\t{\theta}\def\z{\zeta}
\def\P{{\mathbb P}}
\def\Pro{\C\P}					

\def\D{\Delta} \def\G{\Gamma}\def\O{\Omega}
														
\def\ES{\varnothing}  \def\IN{\infty}  
   
  \def\ccal{\mathcal C} \def\ecal{\mathcal E} \def\fcal{\mathcal F}
\def\jcal{\mathcal J}\def\lcal{\mathcal L}
\def\ocal{\mathcal O}\def\rcal{\mathcal R}\def\scal{\mathcal S}
\def\ucal{\mathcal U}\def\wcal{\mathcal W}
   \def\ta{\ccal\rcal}

\def\jn{\!\in\!}  \def\smi{\smallsetminus}

\def\sbs{\!\subset\!}      
\def\sps{\!\supset\!}

\def\Sbs{\!\Subset\!}

\def\Til {\widetilde}
\def \Hat {\widehat}

\def\oli{\overline}				 \def\tms{\times}

\def\benu{\begin{enumerate}} \def\eenu{\end{enumerate}}
\def\beqn{\begin{eqnarray*}}  \def\eeqn{\end{eqnarray*}}
\def\beqnn{\begin{eqnarray}}  \def\eeqnn{\end{eqnarray}}
\def\Sum{\sum\limits_}
\def\debar{\bar\p}

\title[$1$-complete semiholomorphic foliations]{$1$-complete semiholomorphic foliations}

 \author[Samuele Mongodi and Giuseppe Tomassini]{Samuele Mongodi and Giuseppe Tomassini}
 \address{S. Mongodi: Dipartimento di Matematica - Universit\`a di Roma ``Tor Vergata'', Via della Ricerca Scientifica --- I-00133 Roma, Italy}
\email{mongodi@mat.uniroma2.it}
 \address{G. Tomassini: Scuola Normale Superiore, Piazza dei Cavalieri, 7 --- I-56126 Pisa, Italy}
\email{g.tomassini@sns.it}
\date{\today}
\thanks{}

\begin{document}
\maketitle
\begin{abstract}
 A {\it semiholomorphic foliations of type} $(n,d)$ is a differentiable real manifold $X$ of dimension $2n+d$, foliated by complex leaves of complex dimension $n$. In the present work, we introduce an appropriate notion of \emph{pseudoconvexity} (and consequently, $q$-completeness) for such spaces, given by the interplay of the usual pseudoconvexity, along the leaves, and the positivity of the transversal bundle. For $1$-complete real analytic semiholomorphic foliations, we obtain a vanishing theorem for the CR cohomology, which we use to show an extension result for CR functions on Levi flat hypersurfaces and an embedding theorem in $\C^N$. In the compact case, we introduce a notion of weak positivity for the transversal bundle, which allows us to construct a real analytic embedding in $\mathbb{CP}^N$.
\end{abstract}
\renewcommand{\thefootnote}{}
\tableofcontents
\section{Introduction}\nin
A {\it semiholomorphic foliations of type} $(n,d)$, $n\ge1$, $d\ge 1$, is a differentiable real manifold $X$ of dimension $2n+d$, foliated by complex leaves of complex dimension $n$. If $X$ is of class $C^\o$ we say that $X$ is a {\it a real analytic semiholomorphic foliation}. The aim of the present paper is to investigate the geometrical properties of such spaces along the lines of the classical theory of complex spaces. Such spaces were already studied in \cite{GT}, to some extent, but some of the conclusions reached there are spoiled by the lack of a necessary hypothesis in the statements, the positivity of the transversal bundle mentioned below, which allows us not only to correctly state and prove  the results of \cite{GT}, but also to push further our investigation and obtain new insights in the geometry of semiholomorphic foliations; in a similar direction moved also the work of Forstneri\v{c} and Laurent-Thi\'ebaut \cite{fola}, where the attention was focused on the existence of a basis of Stein neighbourhoods for some particular compacts of a Levi-flat hypersurface. We should also mention that, more recently, El Kacimi and Slim\`ene studied the vanishing of the CR cohomology for some particular kinds of semiholomorphic foliation, in \cite{elk, elk1}.
 We are mainly concerned with real analytic semiholomorphic foliations which satisfy some hypothesis of {\em pseudoconvexity}. Pseudoconvexity we have in mind consists of two conditions: $1$-pseudoconvexity along the leaves of $X$, i.e. the existence of a smooth exhaustion function $\phi:X\to\R^+$ which is strongly $1$-plurisubharmonic along the leaves, and positivity of the bundle $N_{tr}$ transversal to the leaves of $X$ (see \ref{trpsfo}). Under these conditions we say that $X$ is $1$-{\it complete}. For $1$-complete real analytic semiholomorphic foliation of type $(n,d)$ we can prove that
 \bit
 \item[1)] every level set $\{\phi\le c\}$ has a basis of Stein neighbourhoods;
 \item[2)] an approximation theorem for CR functions holds on $X$. 
\eit 
 (see Theorem \ref{ST1}).
 
 Once proved these results, the methods of complex analysis apply, in order to study the cohomology of the sheaf $\mathcal{CR}$ of germs of CR functions. We show that, for $1$-complete real analytic semiholomorphic foliation of type $(n,1)$, the cohomology groups $H^q(X,\mathcal{CR})$ vanishes for $q\geq 1$ (see Theorem \ref{ST2}). This implies a vanishing theorem for the sheaf of sections of a CR bundles, which we use in Section 4 to get a tubular neighbourhood theorem and an extension theorem for CR functions on Levi flat hypersurfaces (see Theorems \ref{van}, \ref {tub} and Corollary \ref{ext}). In Section \ref{EMB1} we sketch the proof that a $1$-complete real analytic semiholomorphic foliation of type $(n,d)$ embeds in $\C^{2n+2d+1}$ as a closed submanifold by a smooth CR map (see Theorem \ref{EMB}) and, as an application, we get that the groups $H_j(X;\Z)$ vanish for $j\ge n+d+1$ (see Theorem \ref{omol}). In the last section a notion of {\em weakly positivity} is given for the transversal bundle $N_{tr}$ on a compact real analytic semiholomorphic foliation $X$ of type $(n,1)$ (see \eqref{curv}) and we prove that, under this condition, $X$ embeds in $\C\P^N$ by a real analytic CR map (see Theorem \ref{prem}). 
\section{Preliminaries}\label{pre}
\subsection{$q$-complete semiholomorphic foliations}
We recall that a {\it semiholomorphic foliation of type $(n,d)$} is a (connected) smooth foliation $X$ whose local models are subdomains $U_j=V_j\tms B_j$ of $\C^n\tms\R^d$ and transformations $(z_k,t_k)\mapsto(z_j,t_j)$ are of the form 
\be\label{LT}
\begin{cases}
z_j=f_{jk}(z_k,t_k) \>\> &\\
t_j=g_{jk}(t_k),\>\> &
\end{cases}
\ee
where $f_{jk}$, $g_{jk}$ are smooth and $f_{jk}$ is holomorphic with respect to $z_k$. If we replace $\R^d$ by $\C^d$ and we suppose that $f$ and $h$ are holomorphic we get the notion of {\it holomorphic foliation of type} $(n,d)$. 

Local coordinates $z^1_j,\ldots,z^n_j,t^1_j,\ldots,t^d_j$ satisfying \ref{LT} are called {\it distinguished local coordinates}. 

A continuous function $f:X\to \C$ is a CR function if and only it is holomorphic along the leaves. Let $\ta=\ta_X$ denote the sheaf of germs of smooth CR functions in $X$. 


Given a subset $Y$ of $X$ we denote $\Hat Y_{\ta(X)}$ the envelope of $Y$ with respect to the algebra $\ta(X)$.

A {\it morphism or {\rm CR} map} $F:X\to X'$ of semiholomorphic foliations is a smooth map preserving the leaves and such that the restrictions to leaves are holomorphic. By $\ta(X;X')$ we denote the space of all morphisms $X\to X'$; this space is a closed subspace in $C^{\IN}(X;X')$ so it is a metric space.

A semiholomorphic foliation $X$ is said to be {\it tangentially $q$-complete} if $X$ carries a smooth exhaustion function $\phi:X\to\R^+$ which is $q$-plurisubharmonic along the leaves (i.e. its Levi form has at least $n-q+1$ positive eigenvalues) . In this case we may assume that $\min_X\phi=0$ and $M=\{x\jn X:\phi(x)=0\}$ is a point ($M$ is compact; take $a\jn M$ and a compactly supported nonnegative smooth function $\psi$ such that $\{\psi=0\}\cap M=\{a\}$. If $\e$ is positive and sufficiently small the function $\phi+\e\psi$ has the required properties).

\br\label{mapr}
 Let $X$ be semiholomorphic foliation of type $(n,d)$ which carries a function $\phi$ which is strongly plurisubharmonic along the leaves. Then
 \bit
\item[1)] if $D$ be a bounded domain in $X$, $f:D\to\C$ a CR function continuous up to the boundary ${\rm b}D$ we have
$$
\max\limits_{\oli D}\vert f\vert=\max\limits_{{\rm b} D}\vert f\vert;
$$

\item[2)]  If $f$ is a {\rm CR} function and ${\rm supp}\,f$ is compact then $f=0$.  
The same holds if, more generally, if $\phi$ is $q$-plurisubharmonic along the leaves with $q\le n$.
\eit
\er
\subsection{CR-bundles}\label{crbu}
Let ${\sf G}_{m,s}$ be the group of matrices 
$$
\left( \begin{array}{cc}
	A&B\\
	0&C
\end{array} \right)
$$
where $A\in GL(m;\C)$, $B\in GL(m,s;\C)$ and $C\in GL(s;\R)$.  We set $\sf G_{m,0}=GL(m;\C)$, $\sf G_{0,s}=GL(s;\R)$.
 To each matrix $M\in \sf G_{m,s}$ we associate the linear transformation $\C^m\tms\R^s\rightarrow\C^m\tms\R^s$ given by
 $$
 (z,t)\mapsto (Az+Bt,Ct),
 $$
$(z,t)\in\C^m\tms\R^s$.

Let $X$ be a semiholomorphic foliation of type $(n,d).$  A CR-{\it bundle of type (m,s)} is a vector bundle $\pi:E\to X$ such that the cocycle of $E$ determined by a trivializing distinguished covering $\{U_j\}_j$ is a  smooth CR map $\g_{jk}:U_j\cap U_k\to\sf G_{m,k}$
\be\label{CO}
\g_{jk}=\left( \begin{array}{cc}
	A_{jk}&B_{jk}\\
	0&C_{jk}
\end{array} \right)
\ee
where $C_{jk}=C_{jk}(t)$ is a matrix with smooth entries and $A_{jk}$, $B_{jk}$ are matrices with smooth CR entries. Thus $E$ is foliated by complex leaves of dimension $m+n$ and real codimension $d+s$. If $z_j,t_j$ are distinguished coordinates on $X$ and $\z_j,u_j$ are coordinates on $\C^m\tms\R^s$, then $z_j,\z_j,t_j,\theta_j$ are distinguished coordinates on $E$ with transformations
$$
\begin{cases}
z_j=f_{jk}(z_k,t_k) \>\> &\\
\z_j=A_{jk}(z_k,t_k)\z_k+B_{jk}(z_k,t_k)\theta_k\\
t_j=g_{jk}(t_k)\>\> &\\
u_j=C_{jk}(t_k)\theta_k.
\end{cases}
$$
In particular, $\pi:E\to X$ is a smooth CR map.

The {\it inverse} $E^{-1}$ of $E$ is the CR-bundle $E$ of type $(m,s)$ whose cocycle is  
\be\label{CO1}
\g_{jk}^{-1}=\left( \begin{array}{cc}
	A_{jk}^{-1}&-A_{jk}^{-1}B_{jk}C_{jk}^{-1}\\
	0&C_{jk}^{-1}
\end{array} \right).
\ee
Let $X$ be a semiholomorphic foliation of type $(n,d)$. Then
\bit
\item the tangent bundle $TX$ of $X$ is a CR-bundle of type $(n,d)$; 
\item the bundle $T_\fcal=T_\fcal^{1,0}X$ of the holomorphic tangent vectors to the leaves of $X$ is a CR-bundle of type $(n,0)$; 
\item the complexified $\C TX$ of the tangent bundle is a complex vector bundle.
\eit
\br
If $X$ is embedded in a complex $Z$, its transverse $TZ/TX$ is not a $\sf G_{m,s}$-bundle in general.
\er
\subsection{Complexification}\label{compl}
A real analytic foliation with complex leaves can be complexified, essentially in a unique way: there exists a holomorphic foliation $\Til X$ of type $(n,d)$ with a closed real analytic CR embedding $X\hookrightarrow\Til X$; in particular, $X$ is  a Levi flat submanifold of $\Til X$ (cfr. \cite [Theorem 5.1] {rea}).

In order to construct $\Til X$ we consider a covering by distinguished domains $\{U_j=V_j\tms B_j\}$ and we complexify each $B_j$ in such a way to obtain domains $\Til U_j$ in $\C^n\tms\C^d$. The domains $\Til U_j$ are patched together by the local holomorphic transformations
\be\label{LT2}
\begin{cases}
z_j=\Til f_{jk}(z_k,\tau_k) \>\> &\\
\tau_j=\Til g_{jk}(\tau_k)\>\> &
\end{cases}
\ee
obtained complexifying the (vector) variable $t_k$ by $\tau_k=t_k+i\theta_k$ in $f_{jk}$ and $g_{jk}$ (cfr. \ref{LT}). 

 Let $z_j$, $\tau_j$ distinguished holomorphic coordinates on $\Til U_j$ and $\t_j={\sf Im}\,\tau_j$. Then $\t^s_j={\sf Im}\,\Til g_{jk}(\tau_k)$, $1\le s\le d$, on $\Til U_j\cap\Til U_k$ and consequently, since ${\sf Im}\,\Til g_{jk}=0$ on $X$, 
$$
\t^r=\sum\limits_{s=1}^d\psi^{rs}_{jk}\t^s_k
$$
where $\psi_{jk}=\big(\psi_{jk}^{rs}\big)$ is an invertible $d\tms d$ matrix whose entries are real analytic functions on $\Til U_j\cap\Til U_k$. Moreover, since $\Til g_{jk}$ is holomorphic and $\Til g_{jk|X}=g_{jk}$ is real, we also have $\psi_{jk|X}=\p g_{jk}/\p t_k.$

$\{\psi_{jk}\}$ is a cocycle of a CR-bundle of type $(0,d)$ which extends $N_{\rm tr}$ on a neighbourhood of $X$ in $\Til X$ .
 
In \cite[Theorem 2]{GT} is proved that if $X$ is tangentially $q$-complete then every sublevel $\oli X_c=\{x\in X:\phi\le c\}$ has a fundamental system of neighbourhoods in $\Til X$ which are $(q+1)$-complete complex manifolds.

Let $\Til X$ the complexification of $X$. Then the cocycle of the (holomorphic) transverse bundle $\Til N_{\rm tr}$ (to the leaves of $\Til X$) is 
\be\label{TRB1}
\frac{\p \Til g_{jk}(\tau_k)}{\p \tau_k}=\frac{\p \tau_j}{\p \tau_k}=\left(\frac{\p \tau^\a_j}{\p \tau^\b_k}\right).
\ee\\
\section{Completeness}\label{pos}
 
\subsection{Tranversally complete foliations}\label{trpsfo}
Let $X$ be a semiholomorphic foliation of type $(n,d)$, $n\ge1$, $d\ge 1$ $N_{\rm tr}$ the transverse bundle to the leaves of $X$. A metric on the fibres of $N_{\rm tr}$ is an assignement of a distinguished covering $\{U_j\}$ of $X$ and for every $j$ a smooth map $\l^0_j$ from $U_j$ to the space of symmetric positive $d\tms d$ matrices such that
$$
\l^0_k=\frac{^t\p g_{jk}}{\p t_k}\l^0_j\frac{\p g_{jk}}{\p t_k}.
$$  
Denoting $\p$ and $\oli\p$ the complex differentiation along the leaves of $X$, the local tangential forms
\be\label{FO1}
\o=2\oli\p\p\log\,\l^0_j-\oli\p\log\,\l^0_j\wedge\p\log\,\l^0_j=\frac{\l^0_j\oli\p\p\l^0_j-2\oli\p\l^0_j\wedge\p\l^0_j}{{\l^0_j}^3}
\ee
\be\label{FO2}
\O\oli\p\p\log\,\l^0_j-\oli\p\log\,\l^0_j\wedge\p\log\,\l^0_j=\frac{2\l^0_j\oli\p\p\l^0_j-3\oli\p\l^0_j\wedge\p\l^0_j}{{\l^0_j}^3}
\ee
actually give global tangential forms in $X$. 

The foliation $X$ is said to be {\it tranversally $q$-complete} ({\it strongly tranversally $q$-complete}) if a metric on the fibres of $N_{\rm tr}$ can be chosen in such a way that the hermitian form associated to $\o$ ($\O$) has at least $n-q+1$ positive eigenvalues.
\br
Assume $d=1$ and that $X$ is transversally 1-complete. Then, due to the fact that the functions  $g_{jk}$ do not depend on $z$, $\o_0=\{\p\oli\p\log\,\l^0_j\}$ and $\eta=\{\p\log\,\l^0_j\}$ are global tangential forms on $X$. Moreover, since $\o_0$ is positive and $\o_0=d\eta$, it gives on each leaf $\o_0$ a K\"ahler metric whose K\"ahler form is exact. In particular no positive dimension compact complex subspace can be present in $X$. 
\er
\bex
Every domain $D\sbs X=\C_z\tms\R_u$ which projects over a bounded domain $D_0\sbs\C_z$ is strongly transversally 1-complete. Indeed, it sufficient to take for $\l$ a function $\mu^{-1}\circ\pi$ where $\mu$ is a positive superharmonic on $D_0$ and $\pi$ is the natural projection $\C_z\tms\R_u\to\C_z$.
\eex
\bex\label{EX}
A real hyperplane $X$ in $\C^2$ (or in $\C^n$) is not transversally 1-complete (but clearly it is  an increasing union of strongly transversally 1-complete domains). Indeed, if $X\sbs\C^2$ is defined by $v=0$, where $z=x+iy, w=u+iv$ are holomorphic coordinates, transverse 1-completeness of $X$ amounts to the existence of a positive smooth function $\l=\l(z,u)$, $(z,u)\in\C^2$, such that 
$$
\l\l_{z\oli z}-2\vert\l_z\vert^2>0.
$$
This implies that 
$$
\left(\l^{-1}\right)_{z\oli z}=\frac{2\vert\l_z\vert^2-\l\l_{z\oli z}}{\l^3}<0
$$
so, for every fixed $u$, the function $\l^{-1}$ is positive and superharmonic on $\C_z$, hence it is constant with respect to $z$: contradiction.
\eex 

\bex Let $F:\C^m\to\C^n$ be a holomorphic function and $Y\subseteq\C^n$ be a real analytic manifold of real dimension $k$, consisting only of regular values for $F$. Then, the set $X=F^{-1}(Y)$ can be given the structure of  a  mixed foliation of type $(m-n,k)$; moreover, as $Y$ consists only of regular values of $F$, locally we have that $X\cong Y\times M$ with $M=F^{-1}(y_0)$ for some $y_0\in Y$.

If we take $F$ to be a polynomial map, $M$ embeds in a projective variety $\widetilde{M}\subseteq\mathbb{CP}^m$, therefore, arguing as in the previous example, we have that $X$ cannot be transversally pseudoconvex, because this would amount to the existence of a bounded from below plurisuperharmonic function on $\widetilde{M}$.

On the other hand, let $\psi:\D\to\C^3$ be a proper holomorphic embedding of the unit disc $\D\sbs\C$ and $F:\C^3\to\C^n$ be such that $F^{-1}(0)=\psi(\D)$. There exists $\epsilon>0$ such that $(t,0,\ldots)\in\C^n$ is a regular value for $F$, with $t\in(-\epsilon,\epsilon)$; fix $Y=\{(t,0,\ldots)\in\C^n\ :\ t\in(-\epsilon,\epsilon)\}$ and consider $X=F^{-1}(Y)$.
Now, $X\cong\D\times(-\epsilon,\epsilon)$ and the function $\lambda=(1-|\psi^{Ð1}(z)|^2)^{-1}$ shows that $X$ is transversally 1-complete.
\eex

\bex
Let $X\sbs\C_z^n\tms\R_u$ be the smooth family of $n$-balls of $\C_z^n$, $z=(z_1,\ldots,z_n)$ , $w=u+iv$, defined by
$$
\begin{cases}v=0,&\\ \vert z-a(u)\vert^2<b(u)^2\end{cases}
$$ 
where $a=a(u)$ is a smooth map $\R\to\C^n$, $b=b(u)$ is smooth from $\R$ to $\R$ and $\vert a(u)\vert,\vert b(u)\vert\to+\IN$ as $\vert u\vert\to+\IN$. $X$ is strongly transversally 1-complete with function 
$$
\l(z,u)=\frac{1}{b(u)^2-\vert z-a(u)\vert^2}
$$
and tangentially $1$-complete with exhaustion function $\phi=\l$.
\eex

We want to prove the following 
\bt\label{FUSY}
Let $X$ be a semiholomorphic foliation of type $(n,d)$.  Assume that $X$ is real analytic and strongly transversally 1-complete. Then there exist an open neighbourhood $U$ of $X$ in the complexification $\Til X$  and a non negative smooth function $u:U\to\R$ with the following properties\bigskip
\bit
\item[i)] $X=\{x\in U:u(x)=0\}$\\
\item[ii)] $u$ is plurisubharmonic in $U$ and strongly plurisubharmonic on $U\smi X$.
\eit\bigskip
If $X$ is tranversally 1-complete then property {\rm ii)} is replaced by the following\bigskip
\bit
\item[iii)] the Levi form of the smooth hypersurfaces $\{u={\rm const}\}$ is positive definite.
\eit
\et
\demo
For the sake of simplicity we assume $n=d=1$. Let $\{\l^0_j\}$ a metric on the fibres of $N_{\rm tr}$. With the notations of \ref{compl} let $\Til N_F$ the tranverse bundle on $\Til X$ whose cocycle $\psi_{jk}=\Til g_{jk}$ is defined by (\ref{TRB1}) and let be $\{\mu_j\}$ a smooth metric on the fibres of $N_F'^{-1}$. On $\Til U_j\cap\Til U_k$ we have $\mu_k=\psi_{jk}^2\mu_j$ and consequently, denoting $\mu^0_j$ the restriction $\mu_{j|X}$, we have $\mu^0_j{\l^0_j}^{-1}=\mu^0_k{\l^0_k}^{-1}$ on $U_j\cap U_k$. Thus $\s^0=\{\mu^0_j{\l^0_j}^{-1}\}$ is a smooth section of $N_F'\otimes N_F'^{-1}$, where $N_F'$ is the dual of $N_{\rm tr}$. Extend $\s_0$ by a smooth section $\s=\{\s_j\}$ of $\Til N_F'\otimes \Til N_F'^{-1}$. Then $\{\l_j={\s_j}^{-1}{\mu_j}\}$ is a new metric on the fibres of $\Til N_F$ whose restriction to $X$ is $\{\l^0_j\}$.

Now consider on $\Til X$ the smooth function $u$ locally defined by $\l_j\t_j^2$ (where $\tau_j=t_j+i\theta_j$); $u$ is non negative and positive outside of $X$. Drop the subscript and compute the Levi form $L(u)$ of $u$. We have
\beqnn\label{LE}
L(u)(\xi,\eta)\!\!\!&=&\!\!\!A\xi\oli\xi+2{\sf Re}\,(B\xi\oli\eta)+C\eta\oli\eta=\\
&&\l_{j,z\oli z}\t^2\xi\oli\xi+2{\sf Re}\left\{(\l_{j,z\oli\tau}\t^2+i\l_{j,z}\t)\xi\oli\eta\right\}+\nonumber\\
&&(\l_{j,\tau\oli\tau}\t^2+i\l_{j,\tau}\t-i\l_{j,\oli\tau}\t+\l_j/2)\eta\oli\eta\nonumber
\eeqnn
and 
\be\label{LE1}
AC-\vert B\vert^2=\t^2\left(\l_{j,z\oli z}-2\vert\l_{j,z}\vert^2\right)+\t^3\r
\ee
where $\r$ is a smooth function. The coefficient of $\t^2$ is nothing but that of the form $\l_j^2\O$ (here we denote $\o$ and $\O$, respectively, the forms \ref{FO1}, \ref{FO2} where $\l^0_j$ is replaced by $\l_j$), so if $X$ is strongly transversally pseudoconvex $L(u)$ is positive definite near each point of $X$ and strictly positive away from $X$. It follows that there exists a neighbourhood $U$ of $X$ such that $u$ is plurisubharmonic on $U$.

Assume now that $X$ is transversally 1-complete. The Levi form $L(u)_{|HT(S)}$ of a smooth hypersurface  $S=\{u={\rm const}\}$. $L(u)_{|HT(S)}$ is essentially determined by the function
\beqnn\label{LE2}
&&\t^4\big\{\l_{j,z\oli z}\vert\l_{j,\oli\tau}\t+i\l_j\vert^2-2{\sf Re}\,(\l_{j,z\oli\tau}\t+i\l_{j,\oli z})(\l_{j,\tau}\t-i\l_j)\l_{j,\oli z}+\\
&&(\l_{j,\tau\oli\tau}\t^2+i\l_{j,\tau}\t-i\l_{j,\oli\tau}\t+\l_j/2)\vert\l_{j,z}\vert^2\big\}=\nonumber\\
&&\t^4(\l_j\l_{j,z\oli z}-(3/2)\vert\l_{j,z}\vert^2)+\t^5\r\nonumber
\eeqnn
where $\r$ is a bounded function. The coefficient of $\t^4$ is nothing but that of the form $\l_j^3\o$ so, if $X$ is transversally pseudoconvex, $L(u)_{|HT(S)}$ is not vanishing near $X$. It follows that there exists an open neighbourhood $U$ of $X$ such that the hypersurfaces $S=\{u={\rm const}\}$ contained in $U\smi X$ are strongly Levi convex.

In the general case $\l=\l_j$ is a matrix, $\l_j=(\l_{j,rs})$, so we consider the function $u=\Sum{r=1}^d\l_{j,rs}\t_r\t_s$. Then
\beqnn\label{qLE}
L(u)(\xi,\eta)\!\!\!&=&\!\!\!\Sum{\a,\b=1}^n A_{\a\oli \b}\xi^\a\oli\xi^\b+2{\sf Re}
\Big(\Sum{\a=1}^n\Sum{r=1}^d B_{\a\oli r}\xi^\a\oli\eta^r\Big)+\\
&&\Sum{r,s=1}^d C_{r\oli s}\eta^r\oli\eta^s\nonumber
\eeqnn
where
\be\label{qLE1}
A_{\a\oli \b}=\Sum{\a,\b=1}^n\Big(\Sum{r,s=1}^d\frac{\p^2\l_{j,rs}}{\p z_\a\p\oli z_\b}\t^r\t^s\Big)\xi^\a\oli\xi^\b
\ee
\be\label{qLE2}
B_{\a\oli r}=\Big(\Sum{\a=1}^n\Sum{h=1}^d\frac{\p^2\l_{j,rs}}{\p z_\a\p\oli\tau_r}\t^r\t^s+i\Sum{s=1}^d\frac{\p\l_{j,rh}}{\p z_\a}\t^h\Big)\xi^\a\oli\eta^r
\ee
\be\label{qLE3}
C_{r\oli s}=\Big(\Sum{h,k=1}^d\frac{\p^2\l_{j,hk}}{\p\tau_r\p\oli\tau_s}\t^h\t^k+i\Sum{h=1}^d\frac{\p\l_{j,hs}}{\p\tau_r}\t^h-i\Sum{h=1}^d\frac{\p\l_{j,rh}}{\p \oli\tau_s}\t^h+\l_{j,rs}/2\Big)\eta^r\oli\eta^s
\ee
\enddemo
\br\label{qLE4}
In view of (\ref{LE}), the Levi form $L(u)$ at a point of $X$ is positive in the transversal direction $\eta$.
\er
\bt\label{FSY}
Let $X$ be a semiholomorphic foliation of type $(n,d)$.  Assume that $X$ is real analytic and strongly 1-complete. Then for every compact subset $K\sbs X$ there exist an open neighbourhood $V$ of $K$ in $\Til X$, a smooth strongly plurisubharmonic function $v:V\to\R^+$ and a constant $\bar c$ such that $K\Sbs\{v<\bar c\}\cap X\Sbs V\cap X$.
\et
\demo
Let $\phi:X\to\R^+$ be an exhaustion function, strongly plurisubharmonic along the leaves and a sublevel $X_{c_0}$ of $\phi$ such that $K\sbs X_{c_0}$. Consider $U\sbs\Til X$, $u:U\to\R^+$ as in Theorem \ref{FUSY} and let $v=au+\Til\phi$ where $\Til\phi:U\to\R^+$ is a smooth extension of $\phi$ to $U$ and $a$ a positive constant. Then, in view of Remark \ref{qLE4}, it is possible choose $a$ in such a way to have $L(v)(x)>0$ for every $x$ in a neighbourhood $V$ of $\oli X_{c_0}$, $c_0<c$; this end the proof.
\enddemo
\subsection{Stein bases and a density theorem}
\bt\label{ST1}
Let $X$ be a semiholomorphic foliation of type $(n,d)$.  Assume that $X$ is real analytic and strongly 1-complete and let $\phi:X\to\R^+$ be an exhaustive smooth function strongly plurisubharmonic along the leaves. Then
\bit
\item[i)] $\oli X_c$ is a Stein compact of $\Til X$ i.e. it has a basis of Stein neighbourhoods in $\Til X$.\medskip
\item[ii)] every smooth CR function on a neighbourhood of $\oli X_c$ can be approximated in the $C^\IN$ topology by smooth CR functions on $X$.
\eit
\et
\demo
In view of Theorem \ref{FUSY} we may suppose that $X=\{u=0\}$  where $u:\Til X\to [0,+\IN)$ is plurisubharmonic and strongly plurisubharmonic on $\Til X\smi X$.

Let $U$ be an open neighbourhood of $\oli X_c$ in $\Til X$. With $K=\oli X_c$ we apply Theorem \ref{FSY}: there exist an open neighbourhood $V\Sbs U$ of $\oli X_c$ in $\Til X$, a smooth strongly plurisubharmonic function $v:V\to\R^+$ and a constant $\bar c$ such that $\oli X_c\Sbs\{v<\bar c\}\cap X\Sbs V\cap X$. It follows that for $\e>0$, sufficiently small, $W=\{v<\bar c\}\cap\{u<\e\}\Sbs V\Sbs U$ is a Stein neighbourhood of $\oli X_c$.   

In order to prove ii) consider a smooth CR function $f$ on a neighbourhood $I$ of $\oli X_{c}$ in $X$, and take $c'>c$ such that $\oli X_{c'}\sbs I$. For every $j\in\N$ define $\oli B_j=\oli X_{c'+j}$ and choose a Stein neighbourhod $U_j$ of $\oli B_j$ such that $\oli B_j$ has a fundamental system of open neighbourhoods $W_j\Sbs U_{j+1}\cap U_{j}$ which are Runge domains in $U_{j+1}$. Since $U_0$ is Stein, the $\mathcal O(U_0)$-envelope of $\oli B_0$ coincide with its the plurisubharmonic envelope (cfr. \cite[Theorem 4.3.4]{HO}) hence it is a compact contained in $X\cap U_0$, $\oli B_0$ being the zero set of the plurisubharmonic function $u.$ Thus we may assume that $\oli B_0$ is $\mathcal O(U_0)$-convex. 

We have to prove that for every fixed $C^{k}$-norm on $\oli B_0$ and $\e>0$ there exists a smooth CR-function $g:X\to\C$ such that $\Vert g-f\Vert_{\oli B_0}^{(k)}\le \e.$ 

In view of the approximation theorem of Freeman (cfr. \cite[Theorem 1.3]{FR}), given $\e>0$ there is $\Til f\in \mathcal O(U_0)$ such that $\Vert\Til f-f\Vert_{\oli B_0}^{(k)}<\e$. Now, define $F_0=\Til f$ and for every $j\ge 1$ take $W_j$ such that $\oli W_j$ is a Runge domain in $U_{j+1}$ and a holomorphic function $F_\in\mathcal O(U_{j-1})$ satisfying
$$
\Vert F_{j+1}-F_{j}\Vert_{\oli W_j}^{(k)}<\e/2^{j+1},\>\>j=0,1,\ldots;
$$
then 
$$
F_0+\Sum{j=0}^{+\IN}\big(F_{j+1}-F_j\big)
$$
defines a $C^{\IN}$ function $g$ on $X$ such that $\Vert g-f\Vert_{\oli B_0}^{(k)}\le 2\e$.
\enddemo
\bc\label{env}
Let $X$ be as in Theorem \ref{ST1} and $\oli X_c$, $c>0$. Then, $\oli X_c$ is $\ta(X)$-convex i.e. $\oli X_c$ coincides with its $\ta(X)$-envelope.
\ec
\demo
Set $K=\oli X_c$. In view of Theorem \ref{FUSY}, $X$ is the zero set of a plurisubharmonic function defined in an open neighbourhood $U$ of $X$ in $\Til X$. Let $U_\nu\sbs U$ be a Stein neighborhood of $K$ and ${\rm Psh}(U_\nu)$ the space of the plurisubharmonic in $U_\nu$. Then the envelopes $\Hat K_{\ocal(U_\nu)}$ and $\Hat K_{{\rm Psh}(U_\nu)}$ of $K$ coincide (cfr. \cite[Theorem 4.3.4]{HO}) so in our situation we have 
$$
K\sbs\Hat K_{\ta(X)}\sbs\Hat K_{\ocal(U_\nu)}=\Hat K_{{\rm Psh}(U_\nu)}\sbs U_\nu\cap X.
$$
We obtain the thesis letting $U_\nu\sbs U$ run in a Stein basis of $\oli X_c$. 
\enddemo
\section{Cohomology}\label{coho}
In this section we want to prove the following two theorems
\bt\label{ST2}
Let $X$ be a semiholomorphic foliation of type $(n,1)$. Assume that $X$ is real analytic and strongly 1-complete. Then 
$$
H^q(X;\ta)=0
$$
for every $q\ge 1$.
\et
\bt\label{ST3}
Let $D$ be a connected Stein manifold and $X\sbs D$ a closed, orientable, smooth Levi flat hypersurface. Then
$$
H^q(X;\ta)=0
$$
for every $q\ge 1$.
\et

We need some preliminary facts. Let $\O$ a domain in $\C^n\tms\R^k$ and for every $t\in\R^k$ and set $\O_t=\{z\jn\C^n:(z,t)\jn\O\}$. Following \cite{angr} we say that $\{\O_t\}_{t\jn\R^k}$ is a {\em regular family} of domains of holomorphy if the following conditions are fulfilled:
\bit
\item[a)] $\O_t$ is a domain of holomorphy for all $t\jn\R^k$;
\item[b)] for every $t_0\jn\R^k$ there exist $I_0=\{\vert t-t_0\vert,\e\}$ and a domain $U\sbs\C^n$ such that $\O_{t_0}$ is Runge in $U$ and 
$\cup_{t\jn I_0}\O_t\sbs I_0\tms U$.
\eit
In this situation it can be proved that the sheaf $\ta$ of $\O$ is cohomologically trivial (cfr. \cite[Corollaire pag. 213]{angr}.
\bl\label{xc}
Let $X$ be a $1$-complete semiholomorphic foliation of type $(n,d)$, $\phi:X\to[0,+\IN)$ a  strongly plurisubharmonic along the leaves and $c>0$ a regular value of $\phi$. Then there is $\e>0$ such that the homomorphism
$$
H^q(X_{c+\e};\ta)\rightarrow H^q(X_{c};\ta)
$$
is onto for $q\ge 1$. In particular, the homorphism
$$
H^q(\oli X_{c};\ta)\rightarrow H^q(X_{c};\ta)
$$
is onto for $q\ge 1$. 
\el
\demo
Since $\phi$ is  a strongly plurisubharmonic along the leaves the ``bumps lemma'' method applies: given $x_0\jn{\rm b}X_c$ there exist arbitrarily small neighbourhoods $U\ni x_0$ and open sets $B\Sbs X$ with the following properties:
\bit
\item[i)] $B\sps X_c$, $B\smi X_c\Sbs U$ and $B$ is given by $\{\phi'<0\}$ with $\phi'$ smooth on a neighbourhood of $\oli B$ and strictly plurisubharmonic along the leaves;
\item[ii)] $V:=U\cap B$ and $V\cap X_c$ are CR isomorphic to regular families of domains of holomorphy of $\C^n$;
\eit
Write $B=X_c\cup V$. Since $H^q(V;\ta)=H^q(V\cap X_c;\ta)=0$ for $q\ge 1$, from the Mayer-Vietoris sequence 
$$
\cdots\rightarrow H^q(X_c\cup V;\ta)\rightarrow H^q(X_c;\ta)\oplus H^q(V;\ta)\rightarrow 
$$
$$
H^q(X_c\cap V;\ta)\rightarrow H^{q+1}(X_c\cup V;\ta)\rightarrow\cdots 
$$
we obtain that the homomorphism $H^q(B;\ta)\to H^q(X_{c};\ta)$ is surjective for $q\ge 1$. By iterating this procedure starting from $B$, in a finite number of steps we find a neighbourhood $W$ of $\oli X_c$ such that $H^q(W;\ta)\to H^q(X_{c};\ta)$ is surjective for every $q\ge 1$. It is enogh to take $\e>0$ such that $X_{c+\e}\sbs W$.
\enddemo
{\bf Proof of Theorem \ref{ST2}} We first prove that for every $c$ and $q\ge 1$ the cohomology groups $H^q(X_c;\ta)$ vanish and for this that $H^q(\oli X_c;\ta)=0$ for every $c$ and $q\ge 1$ 
(see Lemma \ref{xc}.)

Let $\Til X$ be the complexification of $X$. By Theorem \ref{ST1} there exists a Stein neighbourhood $U\Sbs\Til X$ of $\oli X_c$ which divides by $X$ into two connected components $U_+$, $_-$ which are Stein domains. Let $\ocal_+$ ($\ocal_-$ be the sheaf of germs of holomorphic functions in $D_+$ ($D_-$) smooth on $U_+\cup X$ $(U_-\cup X)$ and extend it on $U$ by $0$. Then we have on $U$ the exact sequence 
\be\label{ple}
\xymatrix{0\ar[r]&\ocal\ar[r]&\ocal_+\oplus\ocal_-\ar[r]^{\sf re}&\ta\ar[r]& 0}
\ee
where ${\sf re}$ is definined by ${\sf re}(f\oplus g)=f_{|X}-g_{|X}$. Let $\xi$ be a $q$-cocycle with values in $\ta$ defined $\oli X_c$. We may suppose that $\xi$ is defined on $U\cap X$.  Since $U$ is Stein we derive that
$$
H^q(U_+\cup X;\ocal_+)\oplus H^q(U_-\cup X;\ocal_-)\rightarrow H^q(X;\ta)
$$
is an isomorphism for $q\ge 1$. In particular, $\xi=\xi_+-\xi_-$ where $\xi_+$ and $\xi_-$ are represented by two $\oli\p$-closed $(0,q)$-forms $\o_+$ and $\o_-$ respectively which are smooth up to $X$. Now arguing as in \cite[Lemma 2.1]{gase} (replacing $\vert z\vert^2$ by a strictly plurisubharmonic exhaustive function $\phi:U\to\R$) we obtain the following: there exist two bounded pseudoconvex domains $U'_+,\Sbs U$ , $U'_-\Sbs U$ satisfying:
\bit
\item [1)] $U'_+\sbs U_+ $ , $U'_-\sbs U_-$;
\item [2)] $\oli X_c\sbs {\rm b}U'_+\cap X, \oli X_c\sbs {\rm b}U'_-\cap X.$
\eit
By Kohn's theorem (cfr. \cite{kohn}) ${\o_+}_{|U'_+}=\oli\p\a_+$, ${\o_-}_{|U'_-}=\oli\p\a_-$ with $\a_+$, $\a_-$ smoth up to the boundary and this shows that $\xi$ is a $q$-coboundary. Thus $H^q(\oli X_c;\ta)=0$ for every $c$ and $q\ge 1$. 
 
The vanishing of the groups $H^q(X;\ta)$ for $q\ge 2$ now easily follows by a standard procedure. 


The proof of $H^1(X;\ta)=0$ requires approximation of CR functions. Let $X=\bigcup\limits_{s\ge1}^{+\IN}X_s$ and $\Vert \cdot\Vert_{\oli X_s}^{(s)}$ denote the $C^s$ norm, $s\in\N$. Let $\{U_i\}$ be an open covering of $X$ and $\{f_{ij}\}$ an 1-cocycle of $\{U_i\}$ with values in $\ta$. Solve $f_{ij}=f_i^{(1)}-f_j^{(1)}$ on a neighbourhood of $\oli X_1$ and analogously $f_{ij}=g _i^{(2)}-g_j^{(2)}$ on a neighbourhood of $\oli X_2$. Then the $\{g_i^{(2)}-f_i^{(1)}\}$ defines a CR function on a neighbourhood of $\oli X_1$ which can be approximated in the norm $\Vert \cdot\Vert_{\oli X_1}^{(1)}$ by a CR function $h:X\to\C$ (see Theorem \ref{ST1}); thus, defining $f_i^{(2)}=g_i^{(2)}-h$ we have $f_{ij}=f_i^{(2)}-f_j^{(2)}$ on a neighbourhood of $\oli X_2$ and
$$
 \Vert f_i^{(2)}-f_i^{(1)}\Vert_{\oli X_1}^{(1)}<1/2.
$$
Thus, for every $s$ we can find functions $f_i^{(s)}$ solving $f_{ij}=f_i^{(s)}-f_j^{(s)}$ on a neighbourhood of $\oli X_s$ and satisfying 
$$
 \Vert f_i^{(s)}-f_i^{(s-1)}\Vert_{\oli X_s}^{(s-1)}<(1/2)^{s-1}.
$$
This makes the series 
$$
f_i=f_i^{(1)}+\Sum{s=1}^{+\IN}\big(f_i^{(s+1)}-f_i^{(s)}\big)
$$
convergent in $\ta(U_i)$ to a function $f_i$. On $X_m$ we have 
$$
f_i=f_i^{(m)}+\Sum{s=m+1}^{+\IN}\big(f_i^{(s+1)}-f_i^{(s)}\big)
$$
where the sum of the series is a CR function on $X_m$ which is independent of $i$. It follows that$f_i-f_j=f_i^{(m)}-f_j^{(m)}=f_{ij}$ i.e. $\{f_{ij}\}$ is a coboundary and this ends the proof of Theorem \ref {ST2}. 
\enddemo
Observe that the proof of Theorem \ref {ST2} works under the hypothesis of Theorem \ref{ST3} as well.

Clearly, under the hypothesis of Theorem \ref{ST3} if $\{x_\nu\}$ is a discrete subset of $X$ and $\{c_\nu\}$ a sequence of complex numbers   there exist a CR function $f:X\to\C$ such that $f(x_\nu)=c_\nu$, $\nu=1,2,\ldots$. This is actually true under the hypothesis of Theorem \ref{ST2}:
\bt\label{ST4}
Let $X$ be a semiholomorphic foliation of type $(n,1)$. Assume that $X$ is real analytic and strongly 1-complete. Let $A=\{x_\nu\}$ be a discrete set of distinct points of $X$ and $\{c_\nu\}$ a sequence of complex numbers. Then there exists a smooth CR function $f:X\to\C$ such that $f(x_\nu)=c_\nu$, $\nu=1,2,\ldots$. In particular, $X$ is $\ta(X)-convex.$
\et
\demo
Let $I_A$ be the sheaf $\{f\jn\ta:f_{|A}=0\}:$ $\ta/I_A\simeq \prod\C_{x_\nu}$: we have to show that $H^1(X;I_A)=0.$ As usual, we denote $\Til X$ a complexification of $X$, $\phi:X\to\R$ a smooth function displaying the tangential 1-completeness of $X$ and $X_c=\{\phi<c\}$, $c\jn\R.$ In view of Theorem \ref{ST2}we have $H^1(X_;\ta)=0$ for every $c\jn\R$. Moreover, since $A_c:=X_c\cap A$ is finite and $\oli X_c$ has a Stein basis of neighbourhoods in $\Til X$, the map $\ta(X_c)\to\prod_{x\in A_c}\C_x$ is onto; from the exact sequence of cohomology 
$$
\cdots\rightarrow H^0(X_c;\ta)\rightarrow \prod_{x\in A_c}\C_x\rightarrow H^1(X_c;I_A)\rightarrow H^1(X_c;\ta)\to\cdots
$$
it then follows that $H^1(X_c;I_A)=0$. To conclude that $H^1(X;I_A)=0$ we use the same strategy as in Theorem \ref{ST2} taking into account the following: given $x_1,x_2,\ldots,x_k$ in $X_c$ every function $f\jn \ta(X_c)$ vanishing on $x_1,x_2,\ldots,x_k$ can be approximate by functions $f\jn \ta(X)$ vanishing on $x_1,x_2,\ldots,x_k$ (see Theorem \ref{ST1}, ii)).
\enddemo
\bt\label{ST5}
Let $D$ be a connected Stein domain in $\C^n$ (or, more generally a connected Stein manifold) and $X\sbs D$ a closed, orientable, smooth Levi flat hypersurface. Let $\phi:D\to \R$ be exhaustive, strictly plurisubharmonic and $X_c=\{\phi<c\}\cap X$, $c\jn\R$. Then the image of the CR map $\ta(X)\to\ta(X_c)$ is everywhere dense.
\et
\demo
Consider $B$, $U$ as in Lemma \ref{xc}. Observe that, by the bump lemma procedure applied to $\{\phi<c\}$ we may suppose that $B=B'\cap X$ where $B'\Sbs D$ is given by $\{\phi'<0\}$ with $\phi'$ smooth and strictly plurisubharmonic on a neighbourhood of $\oli B'$. In particular, in view of Theorem \ref{ST3}, we have $H^q(B;\ta)=0$ for $q\ge 1$. We are going to prove that the image of $\ta(B)\to\ta(X_c)$ is everywhere dense. Write $B=X_c\cup V$ where $V=B\cap U$. Let $\r:\ta(B)\to\ta(V\cap X_c$: observe that the image of $\r$ is everywhere dense. Let
$$
\s:\ta(X_c)\oplus\ta(V)\to\ta(V\cap X_c)
$$
be the map $f\oplus g\mapsto\r(f)-\r(g)$. In this situation, from the Mayer-Vietoris sequence of Lemma \ref{xc} we deduce that the homomorphism $\s$ is onto. Let $f\jn\ta(X_c)$ and $\r(f)\jn\ta(V\cap X_c)$. Then, since the image of $\r$ is everywhere dense, there exists a sequence $\{g_\nu\}\sbs\ta(V)$ such that $\r(g_\nu)-\r(f)\to 0$. Moreover, in view of Banach theorem applied to $\s$ there are elements $u_\nu\jn\ta(X_c)$, $v_\nu\jn\ta(V)$ with the following properties: $u_\nu\to 0$, $v_\nu\to 0$ and $\r(u_\nu)-\r(v_\nu)=\r(f)-\r(g_\nu)$. Then, setting $f_\nu=g_\nu-v_\nu$ on $V$ and $f-u_\nu$ on $X_c$ we define a sequence of CR functions $f_\nu$ in $B$ such that $f_\nu\to f$. Iterating this procedure starting from $B$ in a finite number of steps we find a domain $W\sps\oli X_c$ with the property that $\ta(W)$ is dense in $\ta(X_c)$ and an $\e>0$ such that $\ta(X_{c+\e})$ is dense in $\ta(X_c)$. Let $I$ the set of $c'>c$ with the property $\ta(X_{c'})$ is everywhere dense in $\ta(X_c)$. Then, an elementary lemma on Fr\'echet spaces shows that $I$ is a closed interval, and moreover, by what is preceding, it is also open. It follows that $\sup I=+\IN$ and (by the same lemma) that $\ta(X)$ is everywhere dense in $\ta(X_c)$.
\enddemo
 A vanishing theorem for compact support cohomology $H^\ast_c(X;\ta$ can be proved in a more general situation. Precisely
 \bt\label{supp}
 Let $X$ be a tangentially 1-complete semiholomorphic foliation of type ${\rm (}n,d{\rm )}$. Then 
 $$
 H^j_c(X;\ta)=0
 $$
 for $j\le n-1$
 \et
 \demo
 Let $\phi:X\to\R^+$ be a smooth function displaying the tangential 1-completeness having only one minimum point $x$ (see \ref{pre}) where $\phi(x)=0$. The bumps lemma method gives the following: if $X_a$ is a sublevel of $\phi$ and $a'>a$ is sufficiently close to $a$ the homomorphism 
 $$
 H^j_c(X_a;\ta)\rightarrow H^j_c(X_{a'};\ta)
 $$
 is an isomorphism for $0\le j\le n-1$. Let $I$ be the set of $a'\ge a$ with such a property. Then $\sup I=+\IN$ so, for every $a$ 
 $$
 H^j_c(X_a;\ta)\rightarrow H^j_c(X;\ta)
 $$
 is an isomorphism for $j\le n-1$. On the other hand, if $a$ is close to $0$ the sublevel $X_a$ is CR isomorphic to a regular family of domains of holomorphy, so $H_c^j(X_a;\ta)=0$ for $0\le j\le n-1$.
 \enddemo
\bc\label{hart}
Let $X$ be a $1$-complete semiholomorphic foliation of type ${\rm (}n,d{\rm )}$ with $n\ge 2$ and $K$ a compact subset such that $X\smi K$ is connected. Then the homorphism
$$
\ta(X)\rightarrow\ta(X\smi K)
$$
is surjective.
\ec 
\subsection{An interpretation}
Let $\Til X$ be the complexification of $X$ and $\wcal$ the sheaf of germs of smooth functions on $\Til X$ which are $\bar\p$-flat on $X$. The restriction map $f\mapsto f_{|X}$ is a surjective homomorphism ${\sf re}:\wcal\to\ta$ (cfr. \cite{nir}). Thus we have an exact sequence

$$
\xymatrix{0\ar[r]&\jcal\ar[r]&\wcal\ar[r]^{\sf re} &\ta\ar[r]& 0}.
$$
We have the following:
\bit
\item[$\a$)] the germs of $\jcal$ are flat on $X$ (cfr. \cite{anhi}) (this is actually true for a generic submanifold of a complex manifold);
\item[$\b$)] the sheaf $\jcal$ is acyclic so for $q  \ge1$
$$
H^q(\Til X;\wcal)\simeq H^q(\Til X;\wcal/\jcal)\simeq H^q( X;\ta). 
$$
\eit
Denote by $\fcal^{(0,q)}$ the sheaf of germs of smooth $(0,q)$-forms in $\Til X$ modulo those which are flat on $X$. Then  
\be\label{FL}
\xymatrix{0\ar[r]& \wcal/\jcal\ar[r]& \fcal^{(0,0)}\ar[r]^{\bar\p}&\cdots\ar[r]^{\bar\p}&\fcal^{(0,n+d)}\ar[r]&0}
\ee
is an acyclic resolution of $\wcal/\jcal$ : given a local $(0,q)$-form $u$ which is $\bar\p$-flat on $X$ there is a local form $(0,q-1)$-form $v$ such that $u-\bar\p v$ is flat on $X$. This can be proved by a direct argument using convolution with Cauchy kernel. Indeed, we can prove something apparently stronger: given a local $\bar\p$-closed $(0,q)$-form $u$, $q\geq 2$, flat on $\{\mathrm{Im}z_n=0\}$, then we can solve the equation $\debar \eta=u$ with a $(0,q-1)$-form $v$, again flat on the same set.

Let $u$ be a local $(0,n)$-form in $\C^n$, flat on $\{\mathrm{Im}z_n=0\}$, then
$$\omega=u_0d\bar{z}_1\wedge \ldots\wedge d\bar{z}_n$$
and we may assume that $u_0$ is compact support. If we set
$$
\eta=u_0\star_1\dfrac{1}{\pi z_1} d\bar{z}_2\wedge \ldots\wedge d\bar{z}_n
$$
(where $\star_1$ denotes the convolution with respect to the variable $z_1$) we obtain a $(0,n-1)$-form, which is again flat on $\{\mathrm{Im}z_n=0\}$, such that $\debar \eta=u$; now, suppose we can solve the problem for $(0,q+1)$-forms by taking convolutions only in the variables $z_1,\ldots, z_{n-1}$, then  given a $(0,q)$-form $u$ in $\C^n$, flat on $\{\mathrm{Im}z_n=0\}$, we consider
$\alpha=u\wedge d\bar{z}_n$: we have
$$\debar\alpha=\debar u\wedge d\bar{z}_n=0$$
and $\alpha$ is flat on $\{\mathrm{Im}z_n=0\}$, so there is a form $\beta=\beta_0\wedge d\bar{z}_n$, flat on the same set, such that
$$\debar \beta_0\wedge d\bar{z}_n=\debar\beta=\alpha=u\wedge d\bar{z}_n\;.$$
We set $\gamma=u-\debar\beta_0$; we have that $\gamma\wedge d\bar{z}_n=0$, so $\gamma=\gamma_0\wedge d\bar{z}_n$, moreover $0=\debar\gamma=\debar\gamma_0\wedge d\bar{z}_n=\debar'\gamma_0\wedge d\bar{z}_n$, if we denote by $\debar'$ the $\debar$ operator in $\C^{n-1}$, with respect to variables $z_1,\ldots, z_{n-1}$.
So, we can solve $\debar'\phi=\gamma_0$, because $q\geq 2$, by taking convolutions with respect to $z_1,\ldots, z_{n-1}$ and treating $z_n$ as a parameter. Therefore, $\phi$ is again flat on $\{\mathrm{Im}z_n=0\}$. In conclusion, we have
$$u=\debar \beta_0 + \debar(\phi\wedge d\bar{z}_n)=\debar \eta\;,$$
with $v$ a $(0,q-1)$-form flat on $\{\mathrm{Im}z_n=0\}$.

Returning to the original problem: given a local $(0,q)$-form $u$ which is $\debar$-flat on $X$, we can solve the equation $\debar\psi=\debar u$ with a flat $(0,q)$-form $\psi$, therefore $\debar(u-\psi)=0$ and consequently $u-\psi=\debar v$ for some $(0,q-1)$-form $v$.

A more general result valid for differential operator with constant coefficients was proved by Nacinovich (cfr. \cite[Proposition 5]{na}). Theorem \ref{ST2} and $\b{\rm)}$ together now imply the following
\bc\label{ST6}
Let $D$ be a connected Stein manifold and $X\sbs D$ a closed, orientable, smooth Levi flat hypersurface. Then
\bit
\item [i)] for every $(0,q)$-form $u$, $q\ge 1$ which is $\bar\p$-flat on $X$ there is a $(0,q-1)$-form $v$ such that $u-\bar\p v$ is flat on $X$;
\item [ii)] for every $(0,q)$-form $u$, $q\ge 2$ which is $\bar\p$-closed and flat on $X$ there is a $(0,q-1)$-form $v$ flat on $X$ and such that $u=\bar\p v.$
\eit 
\ec
\demo
Denote $\ecal^(0,q)(D)$ the space of the $(0,q)$-forms in $D$. In view of $\b$), \ref{FL} and Theorem \ref{ST3}we have
\begin{equation} \label{1} 
H^q(X;\ta)\simeq\frac{\{ u\in \ecal^{(0,q)}(D)\}}{\{\bar\p v+\eta:v\in \ecal^{(0,q-1)}(D),\eta\, {\rm flat\,on\,}X\}}=0;
\end{equation} 
then i) follows.
In order to prove ii) let $u\in \ecal^{(0,q)}(D)$, $q\ge 2$, flat on $X$ and $\bar\p$-closed. Since $D$ is Stein $u=\bar\p \psi$. Then, $\bar\p\psi $ is flat on $X$ therefore, by i), we can write $\psi=\bar\p\eta+v$ where $v$ is flat on $X$. It follows that $u=\bar\p v$ i.e. ii).
\enddemo
\section{Applications}
\subsection{A vanishing theorem for CR-bundles}
\bp\label{fibco}
 Let $X$ be a real analytic semiholomorphic foliation of type ${\rm (}n,1{\rm )}$, $\pi:E\to X$ a real analytic CR-bundle of type ${\rm (}n,0{\rm )}$. Then
\bit
\item[i)] if $X$ is 1-complete $E$ is 1-complete; 
\item [ii)] if $X$ is transversally 1-complete {\rm (}strongly transversally $1$-complete{\rm)} $E$ is transversally $1$-complete {\rm (}strongly transversally $1$-complete{\rm)}.
\eit
\ep
\demo
$E$ is a foliated by complex leaves of dimension $m+n$ and real codimension $1$ (see Section \ref{pre}). If $z_j,t_j$ are distinguished coordinates on $X$ and $\z_j$ is a coordinate on $\C^m$, then $z_j,\z_j,t_j,$ are distinguished coordinates on $E$ with transformations
$$
\begin{cases}
z_j=f_{jk}(z_k,t_k) \>\> &\\
\z_j=A_{jk}(z_k,t_k)\z_k\\
t_j=g_{jk}(t_k).
\end{cases}
$$
Part ii) of the statement is immediate since the transversal bundles of $X$ and $E$ respectively have the same cocycle.
In order to prove part i) let $u:X\to\R$ be a smooth function displaying the tangential $1$-completeness of $X$ and $\g_{jk}:U_j\cap U_k\to \sf G_{m,0}$ the cocycle of $E$ with respect to a trivializing distinguished covering $\{U_j\}_j$ such that $U_j\Sbs X$. Let us consider a hermitian metric $\{h_j\}$ on the fibres of $E$. Then on $U_j\cap U_k$
$$
h_k=^t\!\!\g_{jk}h_j\oli\g_{jk};
$$
denoting $\z_j$ the fibre coordinate on $\pi^{-1}(U_j)$ the family $\{\z_je^uh_j\, ^t\oli\z_j\}_j$ of local functions defines a smooth function $\psi$ on $E$. We drop the subscript $j$ and compute the Levi form $\lcal(\psi)$ of $\psi$ on $\pi^{-1}(U_j)$. One has 
\beqn
\lcal(\psi)(\eta)\!\!\!&=&e^u\Big[\z h\, ^t\oli\z\lcal(u)(\eta) +
\sum_{\a,\b=1}^mh_{\a\oli \b}\,\eta^{\a+n}\,\oli\eta^{\b+n}+\fcal(\eta)\Big]
\eeqn
where $h=(h_{\a\oli\b})$, $\eta=(\eta_1,\ldots,\eta_n,\eta_{n+1},\ldots,\eta_{n+m})$ and $\fcal(\eta)$ is a hermitian form whose coefficients are proportional to $\vert\z\vert$. Arguing as in \cite[Theorem 7]{anna} it is shown that the function $u$ can be chosen in such a way that $\psi$ is plurisubharmonic on $E$ and strictly plurisubharmonic away from the 0-section. Moreover, on the 0-section the Levi form $\lcal(\psi)$ equals the form $\sum_{\a,\b=1}^mh_{\a\oli \b}\,\eta^{\a+n}\,\oli\eta^{\b+n}$ which is strictly positive in the direction of the fibre. It follows that the function $\phi=u\circ\pi+\psi$ is strictly plurisubharmonic everywhere on $E$. This ends the proof of Proposition \ref{fibco}.
\bt\label{van}
Let $X$ be a real analytic semiholomorphic foliation of type ${\rm (}n,1{\rm )}$, $\pi:E\to X$ a real analytic CR bundle of type ${\rm (}m,0{\rm )}$ and $\ecal_{{\small CR}}$ the sheaf of germs of smooth CR-sections of $E$. Assume that $X$ is strongly 1-complete. Then
$$
H^q(X,\ecal_{{\small CR}})=0
$$
for $q\ge 1$.
\et
\demo
Let $p:E^\ast\to X$ be the CR-dual of $E$. If $\g_{jk}:U_j\cap U_k\to \sf G_{m,0}$ is the cocycle of $E$ with respect to a distinguished trivializing open covering $\ucal=\{U_j\}_j$ of $X$, $^t\g_{jk}^{-1}:U_j\cap U_k\to \sf G_{m,0}$ is a cocycle of $E^\star$. By Proposition \ref{fibco} $E^\ast$ is strongly 1-complete. Set $\ta^\ast=\ta_{E^\ast}$ and for each domain $U\sbs X$ let $\ta^\ast(p^{-1}(U))_r$, $r\ge 0$, denote the space of the CR functions on $p^{-1}(U)$ which vanish on $U$ of order at least $r$. We obtain a filtration
 $$
 \ta^\ast(p^{-1}(U))_0\sps\ta^\ast(p^{-1}(U))_1\sps\cdots.
 $$
 If $U$ is a trivializing distinguished coordinate domain and $\z=(\z_1,\ldots,\z_m)$ the fibre coordinate then $f\in\ta^\ast(p^{-1}(U))_r$ if and only if
 $$
 (D^r_\z f)_{|U}:={\frac{\p^rf}{\p\z_1^{r_1}\cdots\p\z_1^{r_m}}}_{|U}=0
 $$
 for $r_1+\cdots+r_m\le r$.
 Moreover, if $U\sbs U_j\cap U_k$ and $\z^j=(\z^j_1,\ldots,\z^j_m)$, $\z^k=(\z^k_1,\ldots,\z^k_m)$ are fibre coordinates on $U_j$, $U_k$ respectively, we have
 $$
 (D^r_{\z^j} f)_{|U}=\g^{(r)}_{jk}(D^r_{\z^k} f)_{|U}
 $$ 
 where $\{{\g_{jk}^{(r)}}\}$ denotes the cocycle of ${\rm S}^{(r)}(E)$, the $r^{\rm th}$ symmetric power of $E$.
 
 Conversely, if $U=\cup_jU_j$ and $\{\s_{j,r_1\cdots r_m}\}$, ${r_1}+\cdots{r_m}=r$, is a smooth CR-section of ${\rm S}^{(r)}(E)$ then the local functions
 $$
 \sum\limits_{r_1+\cdots r_m=r}\s_{j,r_1\cdots r_m}(\z^j_1)^{r_1}\cdots(\z^j_m)^{r_m}
 $$ 
 define a function $f\in\ta^\ast(p^{-1}(U))_{r-1}$.
 
 Thus, for every $r> 0$ we obtain an isomorphism 
 $$
\scal(E)^{(r)}_{{\small CR}}(U)\simeq\ta^\ast(p^{-1}(U))_{r-1}/\ta^\ast(p^{-1}(U))_{r}.
 $$
 where $\scal(E)^{(r)}_{{\small CR}}$ denotes the sheaf of germs of smooth CR-sections of ${\rm S}^{(r)}(E)$.
 Let $\ucal^\ast=\{p^{-1}(U_j)\}_j$: $\ucal^\ast$ is an acyclic covering of $E^\ast$ and from the previous discussion we derive that each cohomology group $H^q(\ucal^\ast,\ta^\ast)$ has a filtration
 $$
 H^q(\ucal^\ast,\ta^\ast)\sps H^q(\ucal^\ast,\ta^\ast)_0\sps H^q(\ucal^\ast,\ta^\ast)_1\sps\cdots
 $$
 with associated graded group isomorphic to
 $$
 G\oplus\bigoplus_{r\ge 0} H^q(\ucal,\scal(E)^{(r)}_{{\small CR}}),
 $$
 for some group $G$. Consequently, by Leray theorem on acyclic coverings we have a filtration on each group $H^q(E^\ast,\ta^\ast)$ with associated graded group isomorphic to
 $$
\bigoplus_{r\ge 0} H^q(X,\scal(E)^{(r)}_{{\small CR}}).
 $$
Since, by Proposition \ref{fibco}, $E^\ast$ is strongly $1$-complete, Theorem \ref{ST2} applies giving $H^q(E^\ast,\ta^\ast)=0$ for every $q\ge 1$, whence $H^q(X,\scal(E)^{(r)}_{{\small CR}})=0$ for every $q\ge 1$, $r\ge 0$. In particular, for $r=0$ we obtain $H^q(X,\ecal_{{\small CR}})=0$ for every $q\ge 1$. 
\enddemo
\subsection{CR tubular neighbourhood theorem and extension of CR functions.}
Let $M$ be a real analytic Levi flat hypersurface in $\C^{n+1}$, $n\ge 1$. In view of \cite[Theorem 5.1]{rea}, there exist a neighbourhood $U\subseteq\C^{n+1}$ of $M$ and a unique holomorphic foliation $\Til\fcal$ on $U$ extending the foliation $\fcal$. A natural problem is the following: given a smooth CR function $f:M\to\C$ to extend it on a neighbourhood $W\sbs U$ by a smoth function $\Til f$ holomorphic along the leaves of $\Til\fcal$. In the sequel we answer this question.

The key point for the proof is the following ``CR tubular neighbourhood theorem'':
\bt\label{tub}
Assume that $M$ is strongly transversally 1-complete. Then there exist an open neighbourhood $W\sbs U$ of $M$ and a smooth map $q:W\to M$ with the properties:
\bit
\item[i)] $q$ is a morphism $\Til\fcal_{|W}\to \fcal$;
\item[ii)] $q_{|M}={\sf id}_M$.
\eit
\et
\demo
Clearly, $M$ is strongly 1-complete. Let $p:N\to M$ the normal bundle of the embedding of $M$ in $\C^{n+k}$. Since $\fcal$ extends on $U$, there is a distinguished open covering $\{U_j\}$ of $U$ with holomorphic coordinates $(z_1^j,\ldots, z_n^j, \tau_1^j)$ such that if $V_j:=U_j\cap M\neq\ES$, then $\{V_j\}$ is a distinguished open covering of $M$ with coordinates $(z^j_1,\ldots, z^j_n, {\sf Re}\tau_1^j)$. In particular, if $\t^j_1$ denotes the imaginary part of $\tau^j_1$, the bundle $N_{|V_j}$ is generated by the vector field 
$$\dfrac{\p}{\p\t^j_1}{\Big|}_{V_j}\;.$$
It is easy to check that $N$ is a $G_{0,1}$-bundle.

We have the following exact sequence of CR-bundles
$$
\xymatrix{0\ar[r]&TM\ar[r]&T\C^{n+k}_{|M}\ar[r]&N\ar[r]& 0}
$$
and, passing to the sheaves of germs of CR morphisms, the exact sequence
$$
\xymatrix{0\ar[r]&\mathcal Hom (N,TM)\ar[r]&\mathcal Hom (N,T\C^{n+k}_{|M})\ar[r]&\mathcal Hom (N, N)\ar[r]& 0}.
$$
Theorem \ref{van} now implies that the homomorphism
$$
\G\big(M,\mathcal Hom (N,T\C^{n+k}_{|M})\big)\longrightarrow\G\big(M,\mathcal Hom (N, N)\big)
$$
is onto.

Let $\phi: N\to T\C^{n+1}_{|M}$ be a CR morphism inducing the identity $ N\to N$. Then, $\phi(\xi)=( p(\xi),\psi(\xi))\in M\tms\C^{n+1}$ where $\xi\to\psi(\xi)$ is a smooth CR map $ N\to\C^{n+1}$ which is of maximal rank along $p^{-1}(x)$ for every $x\in M$.

Then $\xi\mapsto p(\xi)+\psi(\xi)$ defines a smooth CR map $\sigma:N\to\C^{n+1}$, which is locally of maximal rank on $M$,  inducing a smooth CR equivalence from a neighbourhood of $M$ in $N$ and a neighbourhood $W$ of $M$ in $\C^{n+1}$.

We define $q=p\circ\sigma^{-1}$.\enddemo

\bc\label{ext}
Let $M$ be a real analytic Levi flat hypersurface in $\C^{n+1}$, $n\ge 1$, $\fcal$ the Levi foliation on $M$, $\Til\fcal$ the holomorphic foliation extending $\fcal$ on a neighbourood of $M$. Then every smooth CR function $f:M\to\C$ extends to a smooth function $\Til f$ on a neighbourhood of $M$, holomorphic along the leaves of $\Til\fcal$.
\ec
\demo
Take $\Til f=f\circ q$.
\enddemo
\section{An embedding theorem}\label{EMB1}
We want to prove the following
\bt\label{EMB}
Let $X$ be a real analytic semiholomorphic foliation of type $(n,d)$. Assume that $X$ is strongly 1-complete. Then $X$ embeds in $\C^{2n+2d+1}$ as a closed submanifold by a CR map.
\et
First of all let us give the notion of CR polyehedron. Let $X$ be a mixed foliation of type $(n,d)$. A  CR {\em polyhedron of order $N$} of $X$ is an open subset $\sf P\Sbs X$ of the form
$$
{\sf P}=\big\{x\jn X:\vert f_j(x)\vert<1, f_j\jn\ta(X), 1\le j\le N\big\}.
$$
With the notations of Corollary \ref{env} we have the following
\bl\label{pol}
Let $X$ be a real analytic semiholomorphic foliation of type $(n,d)$ strongly 1-complete. Let $\phi:X\to\R^+$ be a smooth function displaying the tangential $1$-completeness of $X$ and $X_c$, $X_{c'}$, $c<c'$ sublevels of $\phi$. Then there exists a CR polyhedron {\sf P} of order $2n+2d+1$ such that $\oli X_c\sbs {\sf P}\sbs X_{c'}$.
\el
\demo
Let $\Til X$ the complexification of $X$. By the proof of Theorem \ref{ST1}, part ii), there exist two Stein domains $U$, $V$ with the following properties: $\oli X_c\sbs U$, $\oli X_{c'}\sbs V$ and $\oli X_c$ has a fundamental system of Stein domains $W\Sbs U\cap V$ which are Runge in $V$. Then take such a $W$ and consider the $\ocal(W)$-envelope $Y$ of $\oli X_c$. By a theorem of Bishop there exists an analytic polyhedron $\sf Q\sbs W$ of order $2n+2d+1$ such that $\oli X_c\sbs\sf Q$ (cfr. \cite[Lemma 5.3.8]{HO}). Since $W$ is Runge in $V$ we can assume that $\sf Q$ is defined by functions $f_1,f_2,\ldots,f_{2n+2d+1}\jn\ocal(V)$. Then, ${\sf P}={\sf Q}\cap X$ is a CR polyhedron of  order $2n+2d+1$ containing $\oli X_c$, defined by CR functions in $V\cap X\sps\oli X_{c'}$: we conclude the proof by approximation (see Theorem \ref{ST1}, ii)).
\enddemo
\nin{\bf Proof of Theorem \ref{EMB}}(Sketch) Let $\phi:X\to\R^+$ be a smooth function displaying the tangential 1-completeness of $X$ and $\Til X$ be a complexification of $X$. Then $X$ is a union of the increasing sequence of domains $X_\nu=\{\phi<\nu\}$ and, by Theorem \ref{ST1}, for every $\nu\jn\N$, $\oli X_\nu$ is a Stein compact and $\ta(X)$ is everywhere dense in $\ta(\oli X_\nu)$. For every let 
$$
F_\nu=\big\{f\jn\ta\big(X;\C^{2n+2d+1}\big):f{\rm\> is\> not\>injective\>and\> regular\> on}\> \oli X_\nu\big\}.
$$ 
Clearly, each $F_\nu$ is a closed subset of the Fr\'echet space $\ta\big(X;\C^{2n+2d+1}\big)$. Moreover, since $\oli X_\nu$ is a Stein compact and $\ta(X)$ is everywhere dense in $\ta(\oli X_\nu)$, $F_\nu$ is a proper subset $\ta\big(X;\C^{2n+2d+1}\big).$ Arguing as in \cite[Lemma 5.3.5]{HO} one proves that no $F_\nu$ has interior points. By Baire's theorem $\cup_{\nu=1}^{+\IN}F_\nu$ is a proper subset of $\ta\big(X;\C^{2n+2d+1}\big)$; in particular there exists $g\jn\ta\big(X;\C^{2n+2d+1}\big)$ which is regular and one-to-one. 

It remains to prove that in $\ta\big(X;\C^{2n+2d+1}\big)$ there exists a map which is regular, one-to-one and proper. 

Following \cite[Theorem 5.3.9]{HO} it is sufficient to construct $f\jn\ta\big(X;\C^{2n+2d+1}\big)$ such that
$$
\big\{x\jn X:\vert f(x)\vert\le k+\vert g(x)\vert\big\}\Sbs X
$$
for every $k\jn\N$. 

The construction of such an $f$ is similar to that given in \cite[Theorem 5.3.9]{HO}, taking into account Remark\ref{mapr}, Theorem \ref{ST1}, Corollary \ref{env} and Lemma \ref{pol}.
\enddemo
\br
The example \ref{EX} shows that the converse is not true, namely a real analytic semiholomorphic foliation embedded in $\C^N$ is not necessarily transversally 1-complete.  
\er
As an application, we get the following
\bt\label{omol}
Let $X$ be a real analytic semiholomorphic foliation of type $(n,d)$.  Assume that $X$ is $1$-complete strongly and transversally 1-complete. Then
$$
H_j(X;\Z)=0
$$
for $j\ge n+d+1$ and $H_{n+d}(X;\Z)$ has no torsion.
\et
\demo
Embed $X$ in $\C^N$, $N={2n+2d+1}$ by a CR map $f$ and choose $a\jn\C^N\smi X$ in such a way that $\psi=\vert f-a\vert^2$ is a Morse function \cite{anfr}. In view of Morse theorem we have to show that no critical points of $\psi$ exists with index larger than $n+d$. The Hessian form $H(\psi)(p)$ of $\psi$ at a point $p\jn X$ is
\beqn
&&\sum_{j=1}^N\Big\vert\sum_{\a=1}^n\frac{\p f_j}{\p z_\a}(p)w_\a\Big\vert^2+\\
&&2\,{\sf Re}\sum_{j=1}^N\Big[\oli f_j(p)-\oli a_j\Big]\sum_{\a,\b=1}^n\frac{\p^2f_j}{\p z_\a\p z_\b}(p)w_\a w_\b+
\sum_{k=1}^dA_k\tau_k
\eeqn
where $A_k$, $1\le k\le d$ is a linear form in $w_1,\ldots,w_n,\oli w_1,\ldots,\oli w_n$. The restriction of $H(\psi)(p)$ to the linear space $\tau_1=\dots=\tau_d=0$ is the sum of a positive form and the real part of a quadratic form in the complex variables $w_1,\ldots,w_n$ so the eigenvalues occur in pairs with opposite sign. It follows that $H(\psi)(p)$ has at most $n+d$ negative eigenvalues.
\enddemo
\br\label{omol1}
In particular, the statement holds for mixed smooth foliations of type $(n,d)$ embedded in some $\C^N$.
\er
\bc\label{lefs}
Let $X\sbs\Pro^N$ be a closed, oriented, semiholomorphic foliation of type $(n,d)$, $V$ a nonsingular algebraic hypersurface which does not contain X. Then the homomorphism
$$
H^j_c(X\smi V;\Z)\rightarrow H^j(X;\Z)
$$
induced by $V\cap X\to X$ is is bijective for $j<n-1$ and injective for $j=n-1$. Moreover, the quotient group
$$
H^{n-1}(V\cap X;\Z)/H^{n-1}(X;\Z)
$$
has no torsion.
\ec
\demo
By Veronese map we may suppose that $V$ is a hyperplane. We have the exact cohomology sequence
\beqn
&&\cdots\rightarrow H^j_c(X\smi V;\Z)\rightarrow H^j(X;\Z\rightarrow
H^j(X\cap V;\Z)\rightarrow\\
&& H^{j+1}_c(X\smi V;\Z)\rightarrow H^{j+1}(X;\Z)\rightarrow\cdots. 
\eeqn
By Poincar\'e duality
$$
H^j_c(X\smi V;\Z)\simeq H_{2n+d-j}(X\smi V;\Z).
$$
Since $X\smi V$ is embedded in $\C^N$, we conclude the proof applying Theorem \ref{omol} and Remark \ref{omol1}. 
\enddemo
\section{The compact case}
Let $X$ be a compact real analytic semiholomorphic foliation of type $(n,1)$. With the notations of \ref{compl} let $\{\psi_{ij}\}$ the cocycle of the CR-bundle of type $(0,1)$ on the complexification $\Til X$ which extends $N_{\rm tr}$. The local smooth functions $h_i$ on $\Til X$ satisfying $h_j=\psi_{ij}^2h_i$ define a metric on the fibres of $N_{\rm tr}$.  We say that  $N_{\rm tr}$ is {\em weakly positive} if a smooth metric $\{h_k\}$ can be chosen in such a way that the $(1,1)$-form
\be\label{curv}
i\p\bar\p \log\phi=i\p\bar\p \log h_k+2i\frac{\p\tau_i\wedge\bar\p\bar\tau_k}{(\tau_i-\bar\tau_k)^2}
\ee
is positive, $i\p\bar\p \log\phi\ge 0$, near $X$, on the complement of $X$.
\bl\label{wcom}
Let $N_{\rm tr}$ be weakly positive and $\phi$ the function on $\Til X$ locally defined by $h_i\psi_{ij}^2$. Then $\phi$ is plurisubharmonic on a neighbourhood of $X$ and its Levi form has one positive eigenvalue in the transversal direction $\tau$. In particular, for $c>0$ small enough, the sublevels $\Til X_c$ of $\phi$ are weakly complete manifolds and give a fundamental system of neighbourhoods of $X$.
\el
\demo
Clearly, $\phi\ge 0$ near $X$ and $\phi>0$ away from $X$. By hypothesis, $i\p\bar\p\log\phi\ge 0$ and
$$
\p\bar\p\phi={\phi}\,\p\bar\p \log \phi+\frac{\p\phi\wedge\bar\p\phi}{\phi};
$$
moreover, locally on $X$ 
$$
2\p\bar\p\phi=h{\p\tau\wedge\bar\p\bar\tau}.
$$
$h=h_k$, $\tau=\tau_k$. The thesis follows since $ h>0$ and $i\p\tau\wedge\bar\p\bar\tau$ is a positive $(1,1)$-form.
\enddemo
\br
Observe that the form \ref{curv} is non negative near $X$ in the complement of $X$, if $X$ admits a space of parameters.
 \er
A CR-bundle $L\to X$ of type $(1,0)$ is said to {\em positive along the leaves} i.e there is a smooth metric $\{h_k\}$ (on the fibres of) $L$ such that\vspace{3mm} 
\be\label{curv1}
\sum\limits_{1\le\a\le n}\frac{\p^2\log h_k}{\p z^k_\a\p\bar z^k_\b}\xi^\a\bar\xi^\b>0.
\ee
\vspace{3mm} 
where $(z^k,t^k)$ are distinguished coordinates.

Under the hypothesis of Lemma \ref{wcom} we have the following  
\bt\label{prem}
Suppose that there exists on $X$ an analytic CR-bundle $L$ of type $(1,0)$ and positive along the leaves. Then $X$ embeds in $\C\P^N$, for some $N$, by a real analytic CR map.
\et
\demo
Extends $L$ on a neighbourhood of $X$ to a holomorphic line bundle $\Til L$ and the metric $\{h_k\}$ to a smooth metric $\{\Til h_k\}$ preserving the condition \ref{curv1} near $X$. For every positive $C\in \R$ consider the new metric $\{\Til h_{k,C}=e^{C\phi}\Til h_k\}$ where $\phi$ is the function defined in Lemma \ref{wcom} and set $\z_1=\z^k_1,\ldots,\z_n=z^k_n, \z_{n+1}=\tau ^k$. 

At a point of $X$ we have\medskip 
\beqn
&&\sum\limits_{1\le\a\le n+1}\frac{\p^2\log \Til h_{k,C}}{\p \z_\a\p\bar \z_\b}\eta^\a\bar\eta^\b=C\!\!\!\!\sum\limits_{1\le\a\le n+1}\frac{\p^2\phi}{\p \z_\a\p\bar \z_\b}\eta^\a\bar\eta^\b+\\&&\sum\limits_{1\le\a\le n+1}\frac{\p^2\log \Til h_k}{\p \z_\a\p\bar \z_\b}\eta^\a\bar\eta^\b=
\lcal_1+\lcal_2.
\eeqn

\vspace{5mm} 
\nin Near $X$ the form $\lcal_2$ is positive for $\vert\eta_{n+1}\vert$ small enough, say $\vert\eta_{n+1}\vert<\e$; in view of Lemma \ref{wcom}, $\lcal_1$ is nonegative and $\lcal_1>0$ for $\vert\eta_{n+1}\vert>0$. It follows that for large enough $C$ the hermitian form $\lcal_1+\lcal_2$ is positive on a neigbourhood of $X$ say $\{\phi<c\}$ for $c$ small enough. 

In this situation a theorem of Hironaka \cite[Theorem 4]{nak} applies to embed $\{\phi<c\}$ in some $\C\P^N$ by a locally closed holomorphic embedding. In particular, $X$ itself embeds in $\C\P^N$ by a CR-embedding. 
\enddemo

\end{document}